\documentstyle[epsf,sectioneqno]{article}
\newtheorem{Algorithm}{Algorithm}[section]
  \title{  EXTENSIONS OF A NEW ALGORITHM FOR THE NUMERICAL SOLUTION OF LINEAR
DIFFERENTIAL SYSTEMS ON AN INFINITE INTERVAL }
\author{    B.M.Brown, M.S.P. Eastham, D.K.R.M$^c$Cormack \\ 
 Department of Computer Science, Cardiff University of Wales, Cardiff, CF2 3XF, U.K. 
  }
\date{}
\begin{document}
\maketitle

\renewcommand{\baselinestretch}{2.0}
\large \normalsize     

 \newcommand{\beq}{\begin{equation}}
\newcommand{\enq}{\end{equation}}
\newcommand{\s}{ {\cal{x}}}
 \newcommand{\Ttildaa}{\overline{T} }
 \newcommand{\Ttilda}{\tilde{T}}

\section{ Introduction}
In two recent papers \cite{BEEM95},\cite{BEM95b}, we have developed an algorithm for the
solution of linear differential systems of the type
\beq
Z'(x)=\rho(x) \{ D + R(x) \} Z(x) \label{eq:1.1}
\enq
on an interval  $[X, \infty)$.      Here $Z$ is an $n$-component vector, $\rho$  is a scalar factor,
$D$ is a constant diagonal matrix and $R$ is a perturbation matrix such that
$R(x)=O(x^{-\delta})$
as $x \rightarrow \infty$,   for some  $\delta > 0$.
     The algorithm implements a repeated transformation
process by means of which  ( \ref{eq:1.1}) is transformed into other systems whose
perturbation matrices are of successively smaller orders of magnitude as $ x \rightarrow \infty$.
When the perturbation reaches a prescribed accuracy, the Levinson
asymptotic theorem \cite{NL48},\cite{MSPE89} provides the solution of the final system in the
process, the solution also possessing that accuracy.  The corresponding
solution of (\ref{eq:1.1}) is then obtained by transforming back.
\par
The algorithm has two main aspects.  The first is the algebraic one of
generating symbolically the matrix terms which appear in the transformation
process.  In \cite{BEEM95} and \cite{BEM95b}, the algorithm was set up in such a way that this
aspect is implemented in the symbolic algebra system Mathematica.   The
other aspect is the numerical one of deriving from the algebra the values
of the solutions of (\ref{eq:1.1}), given  $\rho , D$ and $R$.   Included here are the determination of
explicit error-bounds and procedures for handling large numbers of matrix
terms.

In this paper we develop these ideas in two different but not unrelated
directions. First, we replace the constant matrix $D$ by one which has
entries of differing orders of magnitude as   $x \rightarrow \infty$.   Our methods are
sufficiently indicated by two orders of magnitude, and thus we consider the
system
\beq
Z'(x) = \rho ( x ) \{ \Lambda (x) + R(x) \} Z(x)  \label{eq:1.2}
\enq
where
\beq
\Lambda = \left ( 
\begin{array}{ll}
\lambda \tilde{D} & 0 \\
0 & D 
\end{array} \label{eq:1.3}
\right )
\enq
with constant and diagonal  $\tilde{D} $ and  $D$, while  $\lambda$ is a scalar factor such that $\lambda(x) \rightarrow \infty$ as $ x \rightarrow \infty$.
The amplification of our previous algorithm to cover this new
situation is dealt with in sections 2-4 below. At the end of section 3, we state the Levinson asymptotic theorem in the form that we need, with the necessary accuracy incorporated. As in \cite{BEEM95} and \cite{BEM95b}, a
typical    $\rho(x)$ is $x^\alpha ( \alpha >-1)$.
      However, the value    $\alpha =-1$  also occurs in applications, and
our second extension is to include  $\alpha =-1$  in the scope of the algorithm.  This
matter is the subject of section 5.
\par
Both these extensions occur together in a system (1.2) which is associated
with  the generalised hypergeometric equation.  Since this system can be solved
in terms of special functions, it provides an independent check on the
effectiveness of our algorithm, and this matter is discussed in sections 6 and  7.
\section{ Large scalar factor $\lambda$}
The system (\ref{eq:1.2}) is the first in a repeated transformation process and
therefore we write it as
\beq
Z'_1(x)=\rho(x) \{ \Lambda_1(x) + R_1(x) \} Z_1(x). \label{eq:2.1}
\enq
In order to increase the potential for applications, we broaden  (\ref{eq:1.3}) to
\begin{displaymath}
\Lambda_1(x)=\Lambda_0(x) + \Phi_1(x),
\end{displaymath}
where
\beq
\Lambda_0 = \left (
\begin{array}{ll}
\lambda  \tilde{D} & 0 \\
0 & D 
\end{array}
\right ),\;\;
\Phi_1= \left (
\begin{array}{ll}
  \tilde{\Delta}_1 & 0 \\
0 & \Delta_1 \\
\end{array}
\right ). \label{eq:2.2}
\enq
Here
\begin{displaymath}
\tilde{D}={\rm dg }(d_1,..., d_N), \;\;\; D={\rm dg }(d_{N+1},...,d_n)
\end{displaymath}
are constant matrices, each one having distinct entries,  and
\beq
\lambda(x) \rightarrow \infty \;\;\; ( x \rightarrow \infty).
\enq
Further   $\tilde{\Delta}_1$ and $\Delta_1$  are diagonal with
\beq
\tilde{ \Delta}_1=o(\lambda),\;\;\; \Delta_1 = o(1)\;\;\;(x \rightarrow \infty). \label{eq:2.4}
\enq
The situation considered previously in \cite{BEEM95} and \cite{BEM95b} is where $N = 0$ and, in
extending this earlier work, we shall be brief when there are similarities.
\par
The nature of  $R_1$ is as introduced in \cite[ section 4]{BEEM95}   and  \cite[ section 3]{BEM95b},
that is, $R_1$  is the sum of terms of differing orders of magnitude in the form
\beq
R_1=V_{11}+V_{21}+...+V_{\mu 1} + E_1, \label{eq:2.5}
\enq
where
\beq
V_{11}=o(1), \;\;\;V_{k1}=o(V_{j1})\;\; ( k>j) \label{eq:2.6}
\enq
and
\beq
E_1=o(V_{\mu 1} ) \label{eq:2.7}
\enq
as $ x \rightarrow \infty$. We arrange that
\beq
{\rm dg } V_{11} =0 \label{eq:2.8}
\enq
by transferring any non-zero diagonal terms in  $V_{11} $ to $\Phi_1$.
In (\ref{eq:2.6}) and (\ref{eq:2.7}),  $E_1$  has a prescribed accuracy as  $x \rightarrow \infty$, while the  $V_{j1}$ do not
have that accuracy and they are successively replaced by smaller-order
terms as     we  implement the transformation process which we discuss now.
\par
The sequence of transformations has the form
\beq
Z_m=(I+P_m) Z_{m+1} \;\;\;(m=1,2,...) \label{eq:2.9}
\enq
where  $P_m = o(1)\; ( x \rightarrow \infty )$  and the definition of $P_m$   in (\ref{eq:2.15}) below is an extension of the one
given for $N = 0$ in \cite[(4.11)]{BEEM95}.  Starting with (\ref{eq:2.1}), a typical transformed
system arising from (\ref{eq:2.9}) is
\beq
Z'_m = \rho ( \Lambda_m + R_m ) Z_m \label{eq:2.10}
\enq
where, corresponding to ( \ref{eq:2.5})-( \ref{eq:2.8}), 
\beq
R_m =V_{1m}+V_{2m}+...+V_{\mu m } + E_m \label{eq:2.11}
\enq
with
\begin{displaymath}
V_{km}=o( V_{jm})\;(k >j), \;\;E_m=o(V_{\mu  m})
\end{displaymath}
as $x  \rightarrow \infty $, and dg $V_{1m}=0$.
Here  $\mu$  will depend on $m$, and the process ends at  $m=M$  when (\ref{eq:2.11}) reduces to
$R_M=E_M$
with  $R_M$  thus achieving the prescribed accuracy.
\par
In ( \ref{eq:2.10}) we write
\beq
\Lambda_m = \Lambda_0 + \Phi_m, \;\;\; \Phi_m = {\rm dg } ( \tilde{\Delta}_m , \Delta_m )
\label{eq:2.12}
\enq
where, as we shall see,  $\tilde{\Delta}_m$  and $\Delta_m$ are obtained respectively   from 
 $\tilde{\Delta}_1$  and $\Delta_1$    by adding to them
terms of smaller orders of magnitude.  Thus
\beq
\tilde{\Delta}_m = \tilde{\Delta}_1 + o(\tilde{\Delta}_1), \;\;\; \Delta_m= \Delta_1 +o( \Delta_1)
\label{eq:2.13}
\enq 
as    $x \rightarrow \infty $.
  The transformation (\ref{eq:2.9}) takes (\ref{eq:2.10}) into
\beq
Z'_{m+1} = \rho [ \Lambda_m +(I+P_m)^{-1}
\{ -\rho^{-1} P'_m + V_{1m} + \Lambda_m P_m -P_m \Lambda_m +(R_m-V_{1m} ) +R_m P_m \} ] Z_{m+1}.  \label{eq:2.14}
\enq
We have to define  $P_m$ so that the dominant matrix   $V_{1m}$  is removed from (\ref{eq:2.14})
and, at the same time,  $P_m$  has a form which can be readily included in an
algorithm.
\par
Dropping certain subscripts for clarity, we define the entries  $p_{ij}$  in $P_m$ in terms of the entries
$v_{ij }$ in $V_{1m}$ by 

\beq
p_{ij}=\left \{
\begin{array}{ll}
v_{ij}/\{ \lambda (d_j-d_i)\}\;\; &( 1 \leq i \leq N, 1 \leq j \leq N) \nonumber \\ 
-v_{ij}/( \lambda d_i)\;\; &( 1 \leq i \leq N, N+1 \leq j \leq n) \nonumber \\ 
v_{ij}/( \lambda d_j) \;\; &( N+1 \leq i \leq n, 1 \leq j \leq N) \nonumber \\ 
v_{ij}/ (d_j-d_i)\;\; &( N+1 \leq i \leq n,  N+1 \leq j \leq n)  
\end{array}
\right .
\label{eq:2.15}
\enq
where    $ i \neq  j$, and  dg $P_m=0$.   Then,  by (\ref{eq:2.2})  and (\ref{eq:2.12}), we have
\beq
V_{1m}+\Lambda_mP_m-P_m \Lambda_m=U_m+\Phi_m P_m -P_m\Phi_m
\label{eq:2.16}
\enq
where the entries $u_{ij}$ in $U_m$ are
\beq
u_{ij} =
\left \{
\begin{array}{ll}
 \lambda^{-1}(d_j/d_i)v_{ij}\; &( 1 \leq i \leq  N,\;\;N+1 \leq j \leq n ) \nonumber \\
\lambda^{-1}(d_i/d_j)v_{ij}\; &( N+1 \leq i \leq  n ,\;\; 1 \leq j \leq  N ) 
\label{eq:2.17}
\end{array}
\right .
\enq
and zero otherwise.  We note that  $U_m$   has a factor   $\lambda^{-1}$, making  
$U_m=o(V_{1m})$  as required.
To assess the size of the terms involving   $\Phi_m$ in ( \ref{eq:2.16}),   we first write
\beq
P_m=\tilde{Q}_m+Q_m, \label{eq:2.18}
\enq
where $\tilde{Q}_m $  contains the entries involving   $\lambda$ in (\ref{eq:2.15}) and  $ Q_m $ contains the remaining
entries   $v_{ij}/(d_j-d_i)$.   Then, by (\ref{eq:2.16}) and (\ref{eq:2.12}),
\begin{displaymath}
\Phi_m P_m - P_m \Phi_m = \Ttildaa_m + \Ttilda_m +T_m,
\end{displaymath}
where
\beq
  \Ttildaa_m =  \left ( 
\begin{array}{ll}
\tilde{\Delta}_m &  \\
 & 0 
\end{array}
\right )  \tilde{Q}_m
-
\tilde{Q}_m
\left ( 
\begin{array}{ll}
\tilde{\Delta}_m &  \\
 & 0 
\end{array}
\right ), \label{eq:2.19}
\enq
\beq
\Ttilda_m =  \left ( 
\begin{array}{ll}
0 &  \\
  & \Delta _m  
\end{array}
\right )  \tilde{Q}_m
-
\tilde{Q}_m
\left ( 
\begin{array}{ll}
0 &  \\
  & \Delta_m  
\end{array}
\right ) \label{eq:2.20}
\enq
and
\beq
T_m=  \left ( 
\begin{array}{ll}
0 &  \\
  & \Delta_m  
\end{array}
\right )  Q_m
-
Q_m
\left ( 
\begin{array}{ll}
0 &  \\
  & \Delta_m  
\end{array}
\right ). \label{eq:2.21}
\enq
Again, it follows from (\ref{eq:2.12}), (\ref{eq:2.13}), (\ref{eq:2.15}) and (\ref{eq:2.4}) that 
 $ \Ttildaa_m, \Ttilda_m  $ and $T_m$ are
all  $o(V_{1m})$.     Then, provided also that
\beq
\rho^{-1}P'_m = o(V_{1m}),
\label{eq:2.22}
\enq
we can proceed as in \cite[(4.7)]{BEEM95} and \cite[(2.13)]{BEM95b} to express  $(I+P_m)^{-1}$   as a geometric
series in (\ref{eq:2.14}) and then collate terms of the same order of magnitude to
obtain (\ref{eq:2.10}) - (\ref{eq:2.11}) with $m+1$ in place of $m$.  If the dominant term is
initially denoted by  $S_{m+1}$,  we arrange that   dg $V_{1,m+1}=0$  by defining
\begin{displaymath}
\Lambda_{m+1} =\Lambda_m + {\rm dg } S_{m+1},\;\;\; V_{1,m+1} = S_{m+1} - {\rm dg } S_{m+1}.
\end{displaymath}
In particular,   $\Lambda_m$ is built up from  $\Lambda_1$ by adding terms of successively smaller
orders of magnitude, as forecast by (\ref{eq:2.13}).
\section{ Orders of magnitude}
The transformation process in section 2 is carried out for $m=1,2,...,M-1$      and, as
in \cite{BEEM95},\cite{BEM95b},  it can be formalised as an algorithm if the various $o -$ estimates
are expressed more precisely in terms of suitable orders of magnitude.
Here we consider the situation  where the orders of magnitude  
involve powers of  $x^{-a}$   with a given constant  $ a(>0)$.  More general schemes are
also possible, as in  \cite[ section 3 and example 5.2]{BEM95b}, but we shall not
elaborate further here.
\par
We suppose, then, that $\mu =M-1$ in ( \ref{eq:2.5}) 
and
\beq
V_{j1}=O(x^{-ja}) \;\;\; ( 1 \leq j \leq M-1)
\label{eq:3.1}
\enq
in (\ref{eq:2.6}), while the prescribed accuracy for $E_1$   in (\ref{eq:2.7}) is stated as
\beq
E_1(x) = O(x^{-Ma}).
\label{eq:3.2}
\enq
In practice, (\ref{eq:3.1}) will be an exact order of magnitude except that we are
now allowing the possibility that some  $V_{j1}$ with that exact order may be absent
from (\ref{eq:2.5}).  We make similar order assumptions concerning $ \rho, \lambda$ and $\Lambda_1$
in ( \ref{eq:2.1}) and (\ref{eq:2.2}):
\beq
1/\rho(x) = O(x^{1-Ka}), \label{eq:3.3}
\enq
\beq
1/\lambda(x) = O(x^{-La}), \label{eq:3.4}
\enq
\beq
\tilde{\Delta}_1(x) = \tilde{D_1} +   O(x^{-a}), \label{eq:3.5}
\enq
\beq
\Delta_1(x) =    O(x^{-a}), \label{eq:3.6}
\enq
where  $K$ and  $L$ are positive integers and   $\tilde{D_1}$  is constant.  There is a
further set of assumptions concerning derivatives to state, necessitated by
the presence of  $P'_m$ in (\ref{eq:2.14})  and the definition of   $P_m$  in (\ref{eq:2.15}).  We assume that the
$O-$estimates (\ref{eq:3.1}) and (\ref{eq:3.3}) - (\ref{eq:3.6}) can be differentiated the relevant 
number of times in the obvious way: that is, each differentiation reduces
by $1$ the power of  $x$ on the right - hand side.  The relevant number in this
context is
\begin{eqnarray}
M-j \;\; {\rm for} \;\; V_{j1}, \; M-1 \;\; {\rm for} \;\;  \lambda^{-1}, \nonumber \\
M-2 \;\; {\rm for} \;\;  \rho^{-1},  \tilde{\Delta}_1, \Delta_1.  \nonumber
\end{eqnarray}
\par
A simple induction argument, as in \cite[ section 5]{BEEM95}, then gives
\beq
V_{jm}(x)=O(x^{-(m+j-1)a})
\label{eq:3.7}
\enq
in (\ref{eq:2.11}), where now
\beq
\mu =M-m
\label{eq:3.8}
\enq
and   $E_m = O(x^{-Ma})$,  being the prescribed accuracy  (\ref{eq:3.2}). In addition, by (\ref{eq:2.15}), (\ref{eq:3.3})
and (\ref{eq:3.7}), we have
\beq
\rho^{-1}P'_m =O( x^{-(m+K)a}),
\label{eq:3.9}
\enq
so that (\ref{eq:2.22}) is satisfied.  Further, since  $\tilde{\Delta}_m$ and  $\Delta_m$  also satisfy (\ref{eq:3.5}) and
(\ref{eq:3.6}), other estimates which follow when (\ref{eq:3.7}) is substituted into (\ref{eq:2.17}) -
(\ref{eq:2.21}) are
\begin{eqnarray*}
U_m = O(x^{-(m+L)a}), \;\;\ && \bar{T}_m = O(x^{-(m+L)a}), \\
\Ttilda_m = O(x^{-(m+L+1)a}),\;\;\; && T_m = O(x^{-(m+1)a}).
\end{eqnarray*}
We note that, when $L=1, U_m,\Ttildaa_m$ and $T_m$         all have the same order of magnitude and
can therefore be combined into a single term $\tilde{U}_m$  say.    Then the right-hand
side of (\ref{eq:2.16}) can be expressed simply as
\beq
\tilde{U}_m+\tilde{T}_m
\label{eq:3.10}
\enq
where $\tilde{U}_m= O(x^{-(m+1)a}),$ and $\tilde{T}_m= O(x^{-(m+2)a}).$
  The order estimates now established are the
basis of the algorithm which formalises the transformation process of
section 2.
\par
Before moving on to the formulation of the algorithm, we end this section by stating the Levinson asymptotic theorem as applied to the final system
\beq
Z'_M= \rho ( \Lambda_M+R_M)Z_M \label{eq:3.11}
\enq
in the process (\ref{eq:2.9})-(\ref{eq:2.10}). With $\Phi_M$ as in (\ref{eq:2.12}), we write 
\begin{displaymath}
\Phi_M = {\rm dg } ( \delta_1,...,\delta_n).
\end{displaymath}
Then, by ( \ref{eq:2.12}) and (\ref{eq:2.2}), the theorem \cite[ Theorem 1.3.1]{MSPE89} states that (\ref{eq:3.11}) has solutions $Z_{Mk}\;(1 \leq k \leq n )$ such that, as $x \rightarrow \infty$,
\beq
Z_{Mk}(x)= \{ e_k + \eta_k(x) \} \exp ( \int_X^x \rho (t) \{ \lambda (t) d_k + \delta_k (t) \} dt ) \;\;(1 \leq k \leq N) \label{eq:3.12}
\enq
and $\lambda(t)d_k$ is replaced by $d_k$ in the integral when $N+1 \leq k \leq n$.  Here $e_k$ is the unit coordinate vector in the $k$-direction and $\eta_k(x) \rightarrow0$ as $x \rightarrow \infty$.  The standard conditions in the theorem are that $\rho R_M$ is $L(X,\infty)$ together with the dichotomy condition, which we can take in the form \cite[(1.3.13)]{MSPE89}.
\beq
 {\rm Re}\;  F(x) {\rm has}\; {\rm constant}\; {\rm sign}\; ( {\rm either}\; \leq 0 
{\rm or}\;  \geq 0) \label{eq:3.13}
\enq
  in $[X,\infty),$   where $F(x)$  is any one  of 
\begin{eqnarray*}
\rho \{ \lambda (d_j-d_k ) + \delta_j - \delta_k \} && ( 1 \leq j,k \leq N) \nonumber \\
\rho (  d_j-d_k + \delta_j - \delta_k ) && (N+ 1 \leq j,k \leq n) \nonumber \\
\rho ( \lambda  d_j-d_k + \delta_j - \delta_k ) && ( 1 \leq j\leq N, N+1 \leq k \leq n).\nonumber 
\end{eqnarray*}
If, in addition,
\beq
\mid \int_X^\infty F(x) dx \mid = \infty \;(j \neq k) \label{eq:3.14}
\enq
then, as in \cite[Lemma 3.1]{BEEM95}, we have
\beq
 \mid
\eta_k(x) \mid \leq ( \int_X^\infty \mid \rho(t) \mid \parallel R_M(t) \parallel dt ) /
( 1-n \int_X^\infty | \rho(t) \mid \parallel R_M (t) \parallel dt )
\label{eq:3.15}
\enq
in $[X,\infty)$, where $\parallel R_M \parallel = {\rm max} \mid r_{ijM} \mid \; ( 1 \leq i,j \leq n)$. Thus,
transforming back to the original system (\ref{eq:2.1}) by means of (\ref{eq:2.9}), we see that (\ref{eq:2.1}) has solutions of the form (\ref{eq:3.12}) but pre-multiplied by
\begin{displaymath}
{\Pi}_{m=1}^{M-1} (I+P_m),
\end{displaymath}
and the accuracy achieved in $R_M$ is reflected in $\eta_k$ by (\ref{eq:3.15}). We also note that an obvious situation where (\ref{eq:3.13}) and (\ref{eq:3.14}) hold is given by real $\rho, \lambda$ and $d_j$ in conjunction with (\ref{eq:3.3}) and (\ref{eq:3.4}).
\section{The algorithm for large $\lambda$}
We are now in a position to use the analysis   in the previous   sections   to develop an algorithm and computer code to evaluate  the solutions of  (\ref{eq:1.1}).
The general strategy that we   adopt is broadly similar to that   in \cite{BEEM95}; i.e. we   use the 
asymptotic analysis   in the previous   sections to develop an algorithm and computer code that  will evaluate the solution set in the interval $[X, \infty), \;\;0 < X$,  and then use a numerical ODE solver to compute
the solution over $[0,X]$ with initial conditions at $X$ provided by the asymptotic algorithm.   As in \cite{BEEM95}, a feature of the method is that $X$ is   small, about $10$ or $20$.  However the absolute error in the solution at $X$ is also small:  in the examples that we report on
 it is less than $10^{-6}$.
\par
We   first discuss the algorithm to compute the asymptotic solution of (\ref{eq:1.1}). We   assume that $M-1$ iterations of the algorithm are to be performed.
As in \cite{BEEM95} the procedure consists of   $3$ distinct parts  which we call Algorithms 4.1-4.3.
Algorithm 4.1 is   devoted  to obtaining a set of recurrence formulae which define $S_j$, $j=1,...,M,$ 
in terms of  quantities  with a lower subscript, the only assumption being that the quantities concerned satisfy non-commutative
multiplication. However, the extra structure that
(\ref{eq:1.3}) imposes  must be taken into account in developing the algorithm. 
Also we   can no longer assume that $V_{j1}\;\;(1<j\leq \mu)$ or $E_1$ are initially  zero
and this extra complexity must be taken into account. Further, the approach to estimating the absolute error has been modified. Instead of, as in Algorithm 6.1 of \cite{BEEM95}, computing a cumulative  total error $E_j$, which is the total error committed after $j$ iterations of the procedure, we denote by $E_j$ the error computed as a consequence of  the $j^{th}$ iteration. The total error is computed   from  (\ref{eq:4.1}).
 \begin{Algorithm}

\newcounter{rem1}
\begin{list}%
{( \alph{rem1} )}{\usecounter{rem1}
\setlength{\rightmargin}{\leftmargin}
\setlength{\rightmargin}{\labelwidth}
\setlength{\leftmargin}{\labelwidth}}
\item
Input $M$ to specify the number of iterations and hence the accuracy (\ref{eq:3.2}).
 \item
Initialise    $D_1$,   $V_{j1} \; ( 1 \leq j \leq \mu) $ and $E_1$ with initial   values that define the problem and also define the magnitude of the error.
 \item
For each 
\beq
U
 \in \{ -\rho^{-1} Q^{'}_m, -\rho^{-1} \tilde{Q^{'}}_m, U_m, T_m,\bar{T}_m,\tilde{T}_m,
V_{jm} \;\;
(2 \leq j \leq \mu ), V_{jm}Q_m,V_{jm}\tilde{Q}_m \;(1 \leq j \leq \mu) \}, \label{eq:set}
\enq
determine  the least positive integer $\nu$
with\begin{displaymath}
P^{\nu+1}_m U = O( E_1).
\end{displaymath}
\item
For $r=0$ to $\nu$, determine the order $\sigma_{m+k} =r\sigma_m+($order of $U)$ of $P_m^r U$.
\item
Update $V_{k,m+1}= V_{k,m+1} + (-1)^rP^r_m U $, $E_m=E_m+(-1)^{\mu+1}A_mP^{\mu+1}U$.
\item
Output $S_{m+1} = V_{1,m+1}$.
 \end{list}
\end{Algorithm}
As described in section 2, $\Lambda_{m+1}$ and $V_{1,m+1}$ are  then obtained   from $S_{m+1}$ by moving 
the diagonal entries of $S_{m+1}$ to $\Lambda_{m}$; the  modified matrix $S_{m+1}$ is then renamed $V_{1,m+1}$.
The extra structure that the problem imposes demands that  more quantities appear in the 
algorithm than in the corresponding algorithm of \cite{BEEM95}.
\par
The second algorithm performs the same function as Algorithm  6.2 of \cite{BEEM95}.
\begin{Algorithm}
The quantities
$S_j$ and $E_j$  defined by Algorithm 4.1 are realised as $ n \times n$ matrices.
\end{Algorithm}
The  final algorithm is concerned with the evaluation of the $\sup$ norm of  the total error $E$.
This algorithm has some important differences from its counterpart in \cite{BEEM95}.
\begin{Algorithm}
The total error $E$ is given by  
\beq
E=\sum_{j=1}^M   {\cal {A}}_j E_j (I+{\cal {P}}_j) \label{eq:4.1}
\enq
where we have written 
\begin{eqnarray*}
 {\cal {A}}_j &=& A_j A_{j-1}...A_1 \\
I+{\cal{P}}_j&=& (I+P_1)(I+P_2)...(I+P_j).
\end{eqnarray*}
\end{Algorithm}
In general  both  $A_j$ and  $(I+P_j)$ are small perturbations of $I$.
However  a naive use of the Cauchy-Schwartz  inequality when evaluating the $\sup$ norm of  $E$
will induce a    factor of  $n^{2j }$ into each term involved in the estimation of the norm.
This cost may be considerably reduced   by    noting that
\begin{displaymath}
 {\cal {A}}_j=(I+{\cal {P}}_j)^{-1}
\end{displaymath}
and   using this to obtain a bound for ${\cal {A}}_j$ as 
 \begin{displaymath}
 \parallel {\cal {A}}_j \parallel \leq 1 + \parallel {\cal{P}}_j \parallel/(1-n
 \parallel {\cal{P}}_j \parallel ).      
 \end{displaymath}

 \section {Scalar factor $\rho(x)=x^{-1}$}
By  (\ref{eq:3.9}), the validity of (\ref{eq:2.22}) depends on having  $K \geq 1$  in (\ref{eq:3.3}).  When  $K
= 0$--that is, when  $1/\rho(x) = O(x)$--  those  terms in (\ref{eq:2.15}) which involve   $\lambda$ continue to
satisfy (\ref{eq:2.22}), but those which do not involve $\lambda$ contribute only  $O(xv_{ij}')$   to the
left-hand side of (\ref{eq:2.22}), and such terms are no smaller than   $O(v_{ij})$ .  It
follows that, in order to eliminate terms of the size of  $V_{1m}$  from (\ref{eq:2.14}), the
definition (\ref{eq:2.16}) must be modified to incorporate also   $\rho^{-1}P_m'$ as far as entries
in the lower right-hand block of   $P_m$  are concerned.
\par
Let us suppose then, for simplicity, first that   $\rho(x)=x^{-1}$  and second that a
typical  $(i,j)$ entry in  $V_{1m}$ is $c_{ij}x^{-ma}$,     where  $c_{ij}$  is a constant, and $i$ and $j$  lie in the
range $[N+1,n]$.  Then the $(i,j)$ entry in
\begin{displaymath}
-\rho^{-1}P_m'+V_{1m}+\Lambda_0P_m-P_m\Lambda_0
\end{displaymath}
is zero if, in place of the last of (\ref{eq:2.15}), we define
\beq
p_{ij}=v_{ij}/(d_j-d_i-ma).
\label{eq:5.1}
\enq
The remainder of (\ref{eq:2.15}), involving  $\lambda$,  is  unchanged as is  $U_m$ in (\ref{eq:2.16}) and
(\ref{eq:2.17}).  The definition (\ref{eq:5.1}) is valid provided that
\begin{displaymath}
ma \neq d_j-d_i
\end{displaymath}
for  $1 \leq m \leq M-1$   and all $i $ and $j$  in $[N+1,n]$.  This condition is certainly satisfied
if, for example,
\beq
d_j-d_i < a.
\label{eq:5.2}
\enq
In  terms of (\ref{eq:2.18})-(\ref{eq:2.21}), the modification (\ref{eq:5.1}) leaves unchanged    $\tilde{Q}_m, \Ttildaa_m$ and $\Ttilda_m$,  while  $Q_m$
has  $ma$  subtracted from the denominators of its entries.  The discussion
concerning orders of magnitude in section 3 continues to apply.
\par
In order to obtain a computational algorithm that reflects the analysis of this case we make simple modifications to Algorithms 4.1 and 4.2 of the previous section. The components of (\ref{eq:set})
which involve $Q^{'}_m$ are  omitted, and  in Algorithm 4.2, the definition of $P_m$ is modified in order to
take ( \ref{eq:5.1}) into account. Algorithm 4.3 is unchanged.

\section{The generalised hypergeometric equation : an example}
The generalised hypergeometric equation
\beq
y^{(n)}(x)-x^\alpha \sum_{j=0}^ma_jx^jy^{(j)}(x)=0 \label{eq:6.1}
\enq
provides a situation where the numerical results obtained from our
algorithm can be confirmed independently by the known analytic properties
of the solutions of (\ref{eq:6.1}).  Here  $\alpha$   and the  $a_j$   are constants with   $\alpha > -n$ and $m \leq n-1$.
  On the one hand, (\ref{eq:6.1}) can be written in the form (\ref{eq:2.1}) with
$\rho(x) = x^{-1}$ \cite{MSPE96}.  The algorithm in section 5 is applicable, and solutions with
the asymptotic form given by the Levinson Theorem can be computed back to $x=0$
to a prescribed accuracy.  On the other hand, the analytic theory
developed in \cite{PW86} gives exact relationships between solutions with known
behaviours for small and large $x$, thereby providing an independent check on
the efficiency of our algorithm.
\par
Here we discuss this matter in terms of the example
\beq
y^{(3)}(x)-x^2  y^{(2)}(x)-xy^{(1)}(x)+y(x)=0 \label{eq:6.2}
\enq
 which is the case        $n=3,m=2,\alpha=0$ of (\ref{eq:6.1}).                  We quote from \cite{PW86} for the
detailed results which form the basis of our discussion.  We first note
that
\beq
y_0(x)=x \label{eq:6.3}
\enq
is an exact solution of (\ref{eq:6.2}) and that two other solutions, analytic for
all $x$, have the form
\beq
y_1(x)=\sum_0^\infty c_r x^{3r}, \;\;\; y_2(x)=x^2\sum_0^\infty d_r x^{3r},  \label{eq:6.4}
\enq
 where  $c_0 =d_0=1$.
\par
To express (\ref{eq:6.2}) in the form (\ref{eq:2.1}), we define the column vector    
$Y=(y,y^{(1)},y^{(2)})^T$  and
write  (\ref{eq:6.2})  in the usual way as $Y'=AY$,
where
\begin{displaymath}
A  = 
\left
[ \begin{array}{ccc}
0&1&0 \\
0&0&1 \\
-1&x&x^2 \\
\end{array}
\right ].
\end{displaymath}

Then, as explained in \cite{MSPE96}, we make a series of transformations which, put
together, are
\begin{eqnarray}
Y & = &
\left
[ \begin{array}{ccc}
1&0&0 \\
0&x^{-1}&0 \\
0&0&x \\
\end{array}
\right ]
\left
[ \begin{array}{ccc}
1&1&1 \\
x^{3}&1&-1 \\
x^3-1&0&2x^{-3} \\
\end{array}
\right ]
\left
[ \begin{array}{ccc}
1&0&0 \\
3x^{-3}&1&0 \\
0&0&1 \\
\end{array}
\right ] Z  \nonumber  \\
& = &
\left
[ \begin{array}{ccc}
1+3x^{-3}&1&1 \\
x^2+3x^{-4}&x^{-1}&-x^{-1}\\
x^4-x&0&2 x^{-2} \\
\end{array}
\right ] Z. \label{eq:6.5}
\end{eqnarray}
This gives the $Z-$system
\beq
Z'=x^{-1}( \Lambda + R)Z \label{eq:6.6}
\enq
with 
\beq
\Lambda = {\rm dg }(X-3,1,-1)
\label{eq:6.7}
\enq
\begin{displaymath}
R  = 
3\left
[ \begin{array}{ccc}
-(X-1)^{-1}&0&2(X^2-1)^{-1} \\
4(X-1)^{-1}+4X^{-1}&0&-(X-1)^{-1}-6X^{-1}(X^2-1)^{-1} \\
0&0&(X+1)^{-1} \\
\end{array}
\right ]
\end{displaymath}

and
where   $X=x^3$.  In order to obtain the form (\ref{eq:2.5}), with $\mu=2$  for example, we write
\beq
R=3X^{-1}C_1 + 3X^{-2} C_2 + E_1
\label{eq:6.8}
\enq
 where
\begin{displaymath}
C_1= \left [
\begin{array}{ccc}
-1 &  0 & 0 \\
8 & 0 & -1 \\
0&0&1 
\end{array}
\right ],\;\;
C_2= \left [
\begin{array}{ccc}
-1 &  0 &2 \\
4 & 0 & -1 \\
0&0&-1
\end{array}
\right ],
\end{displaymath}
\begin{displaymath}
E_1= 3\left [
\begin{array}{ccc}
-X^{-2} ( X-1)^{-1} &  0 &2 X^{-2}(X^2-1)^{-1} \\
4X^{-2}(X-1)^{-1}  & 0 & -X^{-2}(X-1)^{-1}-6 X^{-1}(X^2-1)^{-1}  \\
0&0& X^{-2}(X+1)^{-1}
\end{array}
\right ].
\end{displaymath}
On taking the diagonal terms in  $C_1$  over to $\Lambda$,  we finally obtain (\ref{eq:2.1}) with  
$\rho(x)=x^{-1}$  and
\beq
\Lambda_1 = {\rm dg } ( X-3-3X^{-1},1,-1+3 X^{-1} ). \label{eq:6.9}
\enq
In (\ref{eq:2.2}) we have    $\lambda =X=x^3, \tilde{D}=(1),D={\rm dg }(1,-1),
\tilde{\Delta}_1=(-3-3X^{-1}), \Delta_1= {\rm dg }(0,3X^{-1})$.            Thus (\ref{eq:3.4}) - (\ref{eq:3.6}) are satisfied with   $a=3$
and  $ L=1 $ while, by ( \ref{eq:6.8}), (\ref{eq:3.1}) is also satisfied.  The condition (\ref{eq:5.2}) is
also clearly satisfied.  The algorithm in section 5 is therefore
applicable.
\par
The asymptotic solution of (\ref{eq:6.6}) is given by the Levinson Theorem and, when
the first component of (\ref{eq:6.5}) is taken, we obtain solutions 
 $y_{\infty 1 }$  and  $y_{\infty 2 }$    of
(\ref{eq:6.2}) such that
\beq
y_{\infty 1} \sim { x^{-3}} \exp ( \frac{1}{3}x^3), \;\; y_{\infty 2 }\sim { x^{-1}} \label{eq:6.10}
\enq
as  $x \rightarrow \infty $.  These solutions arise from the entries  $X-3$ and      $-1\;$ in $\Lambda$.  The entry
in $1$ in $\Lambda$ also gives rise to a solution asymptotic to $x$, but this solution is
associated with   $y_0(x)$ in (\ref{eq:6.3}).
\subsection{Computational results at $x=0$}
 In order to  test our algorithm, we shall use it to compute the solution $y_{\infty 2}$  at $x=0$.
We do this since we have an independent check on this value provided by a series expansion. The
analysis needed for this is derived in the next section. In this subsection we show how our algorithms 
perform for this example.
\par
We  shall use the theory in the above sections to obtain  the  solution of  (\ref{eq:6.2}) which is asymptotic to $x^{-1}$. The initial conditions required by part (b) of Algorithm 4.1  are obtained from (\ref{eq:6.8}) and are:--
\begin{displaymath}
\mu=2, \;\;V_{11}= 3 X^{-1} ( C_1- {\rm dg} C_1),  \;\;V_{21}= 3X^{-2}C_2,\;\;E_1.
\end{displaymath}
$\Lambda_1$ is determined by (\ref{eq:6.9}).
Algorithm 4.1 gives
\begin{eqnarray}
S_1 &=& V_{11}  \nonumber \\
S_2 &=& - {{{\tilde{Q}}^\prime_1}\over{\rho}} + V_{11} Q_1 + T_1 + \bar{T_1} + U_1 + V_{21}  \nonumber 
\end{eqnarray}
while $S_1$ is obtained directly from (\ref{eq:6.8}).
Algorithm 4.2 realises $S_1$ and $S_2$ as the $3 \times 3$ matrices
 \begin{displaymath}
S_1 = \left ( \begin{array}{ccc}
\frac{-3}{x^3} &  0 & 0\\
 \frac{24}{x^3} & 0 & \frac{-3}{x^3} \\
 0 & 0 & \frac{3}{x^3} \\
\end{array} 
\right )
\end{displaymath}
\beq
S_2 = \left ( \begin{array}{ccc}
\frac{-3}{x^6} &  0 & \frac{6}{x^6}\\
 \frac{36(2 + 7 x^3)}{x^9} & 0 & \frac{-24}{5 x^6} \\
 0 & 0 & \frac{-3}{x^6} \\
\end{array} 
\right ) \label{eq:6.10a}
\enq

The formulae (\ref{eq:2.15}) together with the modification (\ref{eq:5.1}) yield explicit expressions for $P_1$ and $P_2$ respectively.
\par
The Levinson asymptotic theorem ( \ref{eq:3.12}) is applied to the equation (\ref{eq:2.10}) with $m=M=3$ and the result evaluated at $x=10$.
By (\ref{eq:3.2}) therefore, our accuracy at $\infty$ is $O(x^{-9})$. Since we are interested in the solution asymptotic to $x^{-1}$, ( \ref{eq:3.12}) 
( but with $k=n=3$), ( \ref{eq:6.9}) and (\ref{eq:6.10a}) give
\begin{displaymath}
Z_{33}(x) = C \{ e_3 + \eta_3(x) \} \exp 
\left (
\int_X^x t^{-1} ( -1 + 3 t^{-3} -3t^{-6}
\right ).
\end{displaymath}
where  $\eta_3(x)=O(x^{-9})$ by (\ref{eq:3.15}). Here, the arbitrary constant $C$ is introduced and chosen   so that $Z_{33}$ is asymptotic to $x^{-1}$ to within $O(x^{-9})$.
In short
\begin{displaymath}
Z_{33}(x) =  \{ e_3 + \eta_3(x) \}x^{-1} \exp 
\left (
  -x^{-3} +\frac{1}{2}x^{-6}
\right )
\end{displaymath}
and
\beq
Z_{33}(10)= (0 ,0 ,0.0999000999 )^T \label{eq:6.11}
\enq
 while (\ref{eq:4.1})
in this case yields an error estimate
\begin{displaymath}
E(10)=2.09830422 \times 10^{-8}.
\end{displaymath}
Pre-multiplying (\ref{eq:6.11})   by  $(I+P_1)(I+P_2)$ together with (\ref{eq:6.5}) gives  a value
\beq
Y(10)=(0.0999600993,-0.009984070,0.0019920140)^T.\label{eq:6.12}
\enq  
This is then used as the initial conditions to the NAG \cite{Nag:FFT} library numerical ODE solver D02NMF,  which is used to evaluate      
\beq
Y(0)=(1.87778537,-1.76303921,1.99999920)^T \label{eq:6.13}
\enq
which is accurate to an error of  $10^{-7}$  when compared with the result
(\ref{eq:7.3}) from the next section. 
\par
We remark that in order to further verify the accuracy of our method we have implemented the numerical procedure described above in the interval arithmetic code CXSC. This together with the Taylor based ODE solver \cite{Lohphd} has been used to compute enclosures for $Y(0)$ with initial values provided by the asymptotic analysis. These enclosures are
\beq
Y(0) = ( 1.8777^{99}_{72},-1.7630_{30}^{49}, \; ^{2.000011}_{1.999988})
\label{eq:6.14}
\enq
and contain both the computed solution (\ref{eq:6.13}) and the value 
(\ref{eq:7.3}) given by the series in the next section.  We further remark that although the error in the initial condition (\ref{eq:6.12}) at $x=10$ is of order $10^{-8}$, both integrators have 
lost an order of magnitude  in accuracy in integrating over $(10,0)$.
\section{Analytic theory}
We turn now to the analytic theory of (\ref{eq:6.2}) as given in \cite[pp 78-97]{PW86}.  In
the notation of \cite[p.78]{PW86}, we have $n=3,\; p=2, \;\beta_1=-\beta_2=1,\;K=\frac{1}{3}$   and $ \theta= - 1$ in the case of
(\ref{eq:6.2}).  We consider first the solution   $V_{3,2}(x)$  given by \cite[(3.6.4)]{PW86}.  In order
to relate this solution to the real-valued solutions (\ref{eq:6.3}) and (\ref{eq:6.4}), we
define
\begin{eqnarray*}
V(x) &=& \frac{1}{2 \pi} \{ V_{3,2}(x) + V_{3,2} ( x e^{2 \pi i/3}) \} \\
& =&\pi \sum_{k=0}^\infty (-1)^k \frac{ ( 3^{2/3}x)^k \cos( \pi k / 3 )}
{ k ! \Gamma ( \frac{2}{3}-\frac{k}{3})  \Gamma ( \frac{4}{3}-\frac{k}{3}) }.
\end{eqnarray*}
Because of the $\Gamma-$ functions in the denominator, $V(x)$ is a power series in $x^3$     with constant leading term  $\frac{1}{2}3^{3/2}$ ,
together with a single term in  $x$.   Hence, by (\ref{eq:6.3}) and (\ref{eq:6.4}),
\begin{displaymath}
V(x)=\frac{1}{2}3^{3/2}y_1(x)-\frac{\pi 3^{2/3}}{2 \Gamma (\frac{1}{3})}y_0(x).
\end{displaymath}
The asymptotic form of  $V(x)$ is determined by \cite[(3.6.5) - (3.6.6)]{PW86}, and it
is sufficient for our purposes to note that
\begin{displaymath} 
V(x) \sim -\frac{1}{2}3^{3/2}x^{-3} \exp( \frac{1}{3} x^3).
\end{displaymath}
On comparing with (\ref{eq:6.10}) and noting that $y_{\infty 2}$  and  $y_0$ are sub-dominant at $\infty$ , we
conclude that
\beq
y_1(x) =-y_{\infty 1 }(x) + c_1 y_{\infty 2 } (x) + d_1 y_0(x), \label{eq:8.1}
\enq
where  $c_1$  and $d_1$  are constants which may or may not be zero.
\par
Next we consider the solution    $W_{3,2}\;(1,z)$ given by \cite[(3.7.4)]{PW86}.
Again, to keep to real-valued solutions, we define
\begin{displaymath}
W(x)=\frac{1}{\pi}W_{3,2}(1,x e^{2 \pi i /3}) =
\sum_{k=0}^\infty (-1)^k \frac{ \Gamma( \frac{k}{3} + \frac{1}{3})}{k ! \Gamma ( \frac{4}{3}-\frac{k}{3})} ( 3 ^{2/3} x )^k.
\end{displaymath}
As in the case of   $V(x)$,   the $\Gamma-$function  in  the denominator and the
definitions (\ref{eq:6.3}) - (\ref{eq:6.4}), show that
\beq
W(x)=3y_1(x) + \frac{ 3 ^{4/3}}{2 \Gamma ( \frac{2}{3})} y_2(x)-3^{2/3} \Gamma( \frac{2}{3})y_0(x). \label{eq:8.2}
\enq
The asymptotic form of  $ W(x)$ is determined by \cite[ (3.7.5)]{PW86} and again it is
sufficient for our purposes to note that
\begin{displaymath}
W(x) \sim \frac{ 3 ^{1/3}}{\Gamma( \frac{5}{3})}x^{-1}.
\end{displaymath}
It therefore follows from (\ref{eq:6.3}) and (\ref{eq:6.10}) that
\begin{displaymath}
W(x) = \frac{ 3^{1/3}}{\Gamma( \frac{5}{3}) }y_{\infty 2 }(x).
\end{displaymath}
Then  (\ref{eq:8.2}) gives
\begin{displaymath}
y_{\infty 2}(x)=2.3^{-1/3} \Gamma( \frac{2}{3})y_1(x)+y_2(x)-2.3^{-2/3} \{ \Gamma(\frac{2}{3} \}^2y_0(x).
\end{displaymath}
In particular, we have
\beq
y_{\infty 2 }(0)=2.3^{-1/3} \Gamma( \frac{2}{3})= 1.87778588, \label{eq:7.3}
\enq
in agreement with the first component of (\ref{eq:6.13}) and (\ref{eq:6.14}) which is the value of $y_{\infty 2 }(0)$ produced by our algorithm.

\bibliographystyle{plain}
\bibliography{../../bib_dir/bibliography,../../bib_dir/bib1,../../bib_dir/help,../../bib_dir/manual}
 
\end{document}